\documentclass[reqno]{amsart}

\usepackage{amssymb}
\usepackage{latexsym}
\usepackage{amsmath}
\usepackage{verbatim}
\usepackage{amssymb}
\usepackage{amsthm}
\usepackage{amssymb}
\usepackage{bbm}

\usepackage{graphicx}

\numberwithin{equation}{section}

\newcommand{\norm}[2]{\left\| #1 \right\|_{#2}} 
\newcommand{\abs}[1]{\left|#1\right|} 

\newcommand{\R}{{\mathbb {R}}}

\newtheorem{theorem}{Theorem}
\newtheorem*{theorem*}{Theorem}

\newtheorem{proposition}{Proposition}[section]

\theoremstyle{remark} 
\newtheorem{remark}[proposition]{\bf Remark}

\title[Kinetic model of chemotaxis]
{Global existence for a kinetic model of chemotaxis via dispersion and Strichartz estimates. }

\author[N.Bournaveas]{Nikolaos Bournaveas}

\address{University of Edinburgh, School of Mathematics, JCMB, King's Buildings, Edinburgh EH9 3JZ, UK}

\email{N.Bournaveas@ed.ac.uk}

\author[V.Calvez]{Vincent Calvez}

\address{\'Ecole Normale Sup\'erieure, D\'epartement de math\'ematiques et applications, CNRS UMR8553, 45 rue d'Ulm,
F 75230, Paris, cedex 05, France}

\email{Vincent.Calvez@ens.fr}

\author[S.Guti\'errez]{ Susana Guti\'errez}

\address{University of Nottingham, School of Mathematical Sciences, Division of Theoretical Mechanics, 
University of Nottingham, University Park, Nottingham NG7 2RD, UK} 

\email{susana.gutierrez@nottingham.ac.uk}

\author[B.Perthame]{Beno\^it Perthame}

\address{\'Ecole Normale Sup\'erieure, D\'epartement de math\'ematiques et applications, CNRS UMR8553, 45 rue d'Ulm,
F 75230, Paris, cedex 05, France}

\email{Benoit.Perthame@ens.fr}

\date{\today}

\thanks{  }  


\begin{document}

\begin{abstract} We investigate further the existence of solutions to kinetic models of chemotaxis. These are nonlinear transport-scattering equations with a quadratic nonlinearity which  have been used to describe the motion of bacteria since the 80's when experimental observations have shown they move by a series of  'run and tumble'. The existence of solutions has been obtained in several papers \cite{CMPS, HKS1, HKS3} using direct and strong dispersive effects.

Here, we use the weak dispersion estimates of \cite{CP} to prove  global existence in various situations depending on the turning kernel. In the most difficult cases, where both the velocities before and after tumbling appear, with the known methods, only Strichartz estimates can  give a result, with a smallness assumption. 
\end{abstract}

\maketitle
 
\vspace{1cm}

\section{Introduction}
In this paper we study the Othmer-Dunbar-Alt kinetic model of chemotaxis and prove global existence
of solutions under various assumptions on the turning kernel. This model was proposed in \cite{alt, ODA}
for the description of the chemotactic movement of cells in the presence of a chemical substance and
it can be thought of as the mesoscopic analogue of the famous Keller-Segel model \cite{KS1, KS2, KS3, Horstmann03}. It was proposed in the 80's, after the experimental observation that bacteria ({\em E. Coli} in the present case, but this is also true for other bacteria as {\em B. Subtilis} for instance) move by a series of  'run and tumble' corresponding to the clockwise or counterclockwise activations of their flagellas in response to chemoattractant substances and receptors saturation.

Denoting the cell density by $f(t,x,v)$ and the density of the chemoattractant by $S(t,x)$ the equations
read as follows:
\begin{subequations}\label{Alt}
\begin{align}
&\partial_{t} f + v\cdot\nabla_{x} f = \int_{V}\left( T[S] f' - T^{*}[S] f \right) dv', \label{Alt1}\\
& f(0,x,v)= f_{0}(x,v), \label{Alt2}\\
&\beta S -\Delta S = \rho:=\int_{V} f(t,x,v) dv , \ \beta = 0, 1. \label{Alt3}
\end{align}
\end{subequations}
 We have used the abbreviations 
$\int_{V} T[S] f' dv'= \int_{V} T[S](t,x,v,v') f(t,x,v') dv'$ and 
$\int_{V} T^{*}[S] f  dv' = \int_{V} T[S](t,x,v',v) f(t,x,v) dv'$.  
The velocity space $V$ is assumed to be a bounded three dimensional domain, typical examples being
balls $\{|v|\leq R\}$ and spherical shells $\{r\leq |v| \leq R\}$.
As a consequence, if $f_{0}$ has compact support in $x$, this will be so for all later times, and many aspects
of the present paper can be simplified or improved, but of course to the expense of generality.
Therefore we do not go in that direction.

Several  earlier works have been devoted to the mathematical study of this kinetic model of chemotaxis. In \cite{HO}, the linear system has been studied (i.e. with a given field $S$) and in particular a major issue has been exhibited concerning the 'memory' effect present in the model through a time scale $\epsilon$ in expressions as  $S(x-\epsilon v')$ or $S(x+\epsilon v)$. Not only this is a major experimental observation related to receptors saturation, but it also is responsible for an asymmetric kernel (in $v$, $v'$) which yields the drift term in the Keller-Segel model that is derived in the diffusion limit of equation (\ref{Alt}).  The meaning of $S(t,x-\epsilon v')$ is that cells measure the concentration of the chemical $S$ at position $x-\epsilon v'$ before changing their direction at position $x$, because of an internal memory effect. The other contribution $S(t,x+\epsilon v)$ is interpreted as follows: cells are able to measure the concentration at a location $x+\epsilon v$ thanks to sensorial protrusions. We set $\epsilon = 1$ in the following without loss of generality.  

The nonlinear Initial Value Problem  \eqref{Alt} was first studied in \cite{CMPS} where global existence was proved in $d=3$
dimensions  under the assumption that the turning kernel satisfies the condition
\[0\leq T[S](t,x,v,v') \leq C \Big( 1 + S(t,x+v) + S(t,x-v') \Big) \]
and the initial data satisfy $0\leq f_{0} \in L^{1}(\R^6)\cap L^{\infty}(\R^6)$. 
The proof starts with the fact
that the $L^{1}_{x,v}$-norm of the solution $f$ is a-priori bounded thanks to conservation of mass
\begin{equation}\label{mass}
\iint_{V} f(t,x,v) dv dx = \iint_{V} f_{0}(x,v) dv dx {\bf =: M,}
\end{equation}
and then proceeds to bootstrap higher $L^{p}_{x,v}$-norms based on strong dispersion estimates (see \cite {Gl,P1}). 

The same method was used in  the paper \cite{HKS1} which points out the difference in the dispersive arguments for terms involving both $S(t,x+v)$ and $S(t,x-v')$. The authors
 prove global existence in $d=3$ dimensions under the assumption
\[0\leq T[S](t,x,v,v') \leq C \Big( 1 + S(t,x+v) + |\nabla S(t,x+v)| \Big) \]
or
\[0\leq T[S](t,x,v,v') \leq C \Big( 1 + S(t,x-v') + |\nabla S(t,x-v')| \Big) , \]
and in $d=2$ dimensions together with $\beta=1$ under the assumption
\begin{multline*}
0\leq T[S](t,x,v,v') \leq C \Big( 1 + S(t,x+v) + S(t,x-v')   \\
  + |\nabla S(t,x+v)| + |\nabla S(t,x-v')| \Big) .
\end{multline*}
The main difficulty appears: scattering terms involving $S(x-v')$ or $S(x+v)$ lead to use two different dispersion estimates, that lead to use a bootstrap with integrability exponents that are only compatible in dimensions less than four. The same dispersive method has been pushed forward in \cite{HKS1, HKS3}, including more general biologically relevant turning kernels and pointing out several limitations. 

For more  results and models involving kinetic equations, see \cite{FLP, charod,HillenPS}, for hyperbolic models 
 \cite{HKS2, EH} and for surveys on the kinetic aspects   \cite{P1, P2}.

In this paper we use the dispersion and Strichartz estimates for solutions of the kinetic
transport equation proved in  \cite{CP} to extend the three dimensional results of \cite{CMPS} and 
\cite{HKS1} to more general turning kernels. Compared to \cite{CMPS} and \cite{HKS1} where 
$L^{p}_{x,v}$-spaces are used, the main feature of our present estimates is to work in $L^{p}_{x}L^{q}_{v}$-spaces
for appropriate choices of $p$ and $q$ (see remark \ref{rem:Admissible exponents} for instance).

In our first result we combine the dispersion estimate of \cite{CP} with  the well-known 
consequence of Calder\'on-Zygmund theory that any second derivative can be controlled in $L^p$
($1<p<\infty$) by the Laplacian in $L^p$, to prove global existence for the IVP \eqref{Alt}
under assumption \eqref{hyp2} below. The latter assumption allows the turning kernel $T[S]$ to be controlled by
second derivatives of the chemoattractant density $S$. Notice that this result is valid in all
dimensions $d\geq 2$.

\smallskip

\begin{theorem}\label{thm2}
Let $d\geq 2$ and suppose that the (continuous) turning kernel satisfies 
\begin{equation}\label{hyp2}
0\leq T[S](t,x,v,v') \leq C \left( 1 + \sum_{|\alpha|\leq 2} \abs{\partial^{\alpha} S(t,x+v)} \right) .
\end{equation}
Fix $p\in \left(1,\frac{d}{d-1}\right)$. If the initial
data $f_{0}\in L^{1}(\R^{2d})$ is nonnegative and such that 
$\norm{f_{0}(x-tv,v)}{L^{p}(\R^{d}_{x} ; L^{1}(\R^{d}_{v}))}$ is finite for all $t > 0$ (\footnote{This assumption means simply that the solution 
$f(t,x,v)=f_{0}(x-tv,v)$ of the linear homogeneous kinetic transport equation with initial data 
$f_{0}$ belongs to $L^{p}_{x}L^{1}_{v}$ for all times. This assumption is satisfied if for example 
$f_{0}$ is $L^{1}_{x}L^{p}_{v}$ (see proposition \ref{disp}).}),
then  the Cauchy problem \eqref{Alt} with $\beta=1$ has a  global weak solution $f$ with 
$ f(t)\in L^{1}(\R^{2d}) \cap L^{p}(\R^{d}_{x} ; L^{1}(V)) $
for all $t\geq 0$. 
\end{theorem}

As it was already commented in \cite{CMPS} and \cite{HKS1} it is difficult to mix terms involving $x+v$ with terms
involving $x-v'$. In this direction we shall prove the following result.

\begin{theorem}\label{thm1}
Let $d=3$ and suppose that the (continuous) turning kernel satisfies 
\begin{equation}\label{hyp1}
0\leq T[S](t,x,v,v') \leq C \Big( 1+   S(t,x+v) + S(t,x-v') + \abs{\nabla S (t,x+v)}\Big) .
\end{equation}
Let $q\in ( 1 , 3/2)$. Then there exists a $p\in (3/2,3)$ (depending on $q$) such that if the initial
data $f_{0}\in L^{1}(\R^6)$ is nonnegative and such that 
$\norm{f_{0}(x-tv,v)}{L^{p}(\R^{3}_{x} ; L^{q}(\R^{3}_{v}))}$ is finite for 
all $t >0$,
then  the Cauchy problem \eqref{Alt} has a  global weak solution $f$ with 
$ f(t)\in L^{1}(\R^6) \cap L^{p}(\R^{3}_{x} ; L^{q}(V)) $
for all $t\geq 0$.
\end{theorem}

\smallskip

Hypothesis \eqref{hyp1} does not allow putting together the two gradients 
$\nabla S(t,x+v)$ and $\nabla S(t,x-v')$. However, our next result shows 
that  if we add the assumption that the critical $L^{3/2}_{x,v}$-norm of 
the initial data is sufficiently small then we have global existence 
under a very general hypothesis on the turning kernel, see \eqref{hyp3} 
and the even weaker \eqref{hyp4}. 
The proof uses the Strichartz estimates of \cite{CP}  
and can be made to work in $d=2$ and $4$ dimensions too,  
see Remark \ref{numerology4}.

\begin{theorem}\label{thm3} Let $d=3$. Consider nonnegative initial data 
$f_{0}\in L^{1}(\R^{2d}) \cap L^{a}(\R^{2d})$, where $\frac{3}{2} \leq a \leq 2$,  
and  assume that $\norm{f_0}{L^{a}(\R^{2d})}$ is sufficiently small. 
Assume that the (continuous) turning kernel $T[S]$ satisfies
the condition
\begin{equation}\label{hyp3}
0\leq T[S](t,x,v,v')\lesssim \abs{S(t,x \pm v)} + 
\abs{S(t,x \pm v')}+ \abs{\nabla S(t,x \pm v)} + \abs{\nabla S(t,x \pm v')}
\end{equation}
where any combination of signs is allowed in the right hand side. 
Then the IVP \eqref{Alt} with $\beta=1$
has a global weak solution $f\in L^{3}_{t}\left([0,\infty); L^{p}\left(\R^{3}_{x} ; L^{q}(V)\right)\right) $, 
where $\frac{1}{p}=\frac{1}{a}-\frac{1}{9}$ and $\frac{1}{q}=\frac{1}{a}+\frac{1}{9}$.
This result also holds if hypothesis \eqref{hyp3} is replaced by the weaker \eqref{hyp4} below.
\end{theorem}

\section{Dispersion and Strichartz estimates}\label{DSE}
In this section we collect the dispersion and Strichartz estimates we shall use later.
We start with the dispersion estimate.

\begin{proposition}(Dispersion estimate, \cite{CP}) \label{disp} 
Let $f_{0}\in L^{q}(\R^{d}_{x} ; L^{p} (\R^{d}_{v}))$ where 
$1\leq q \leq p \leq \infty$, and let $f$ solve 
\begin{equation}\label{kte}
\partial_{t} f + v \cdot \nabla_{x} f =0 
\end{equation}
with initial data $f(0,x,v)=f_{0}(x,v)$. Then 
\begin{equation}\label{dispest1}
\norm{f(t)}{L^{p}(\R^{d}_{x} ; L^{q} (\R^{d}_{v}))} 
 \leq \frac{1}{|t|^{d\left(\frac{1}{q}-\frac{1}{p}\right)}}
\norm{f_{0}}{L^{q}(\R^{d}_{x} ; L^{p} (\R^{d}_{v}))} .
\end{equation}
\end{proposition}

We are going to need the following two versions of the dispersion estimate. First of all observe that the solution
 of \eqref{kte} with initial data $f_{0}(x,v)$ is simply $f(t,x,v)=f_{0}(x-tv,v)$. Therefore the dispersion estimate
says that for any function $h\in  L^{q}(\R^{d}_{x} ; L^{p} (\R^{d}_{v}))$, where $1\leq q \leq p \leq \infty$,
we have
\begin{equation}\label{dispest2}
\norm{h(x-tv,v)}{L^{p}(\R^{d}_{x} ; L^{q} (\R^{d}_{v}))}
\leq \frac{1}{|t|^{d\left(\frac{1}{q}-\frac{1}{p}\right)}}
\norm{h(x,v)}{L^{q}(\R^{d}_{x} ; L^{p} (\R^{d}_{v}))} .
\end{equation}
Replacing $h(x,v)$ by $h(x,v)\mathbbm{1}_{V}(v)$ we get
\begin{equation}\label{dispest3}
\norm{h(x-tv,v)}{L^{p}(\R^{d}_{x} ; L^{q} (V))}
\leq \frac{1}{|t|^{d\left(\frac{1}{q}-\frac{1}{p}\right)}}
\norm{h(x,v)}{L^{q}(\R^{d}_{x} ; L^{p} (V))} .
\end{equation}
In the special case of a  function $h(x)$ which is independent of $v$ we get
\begin{equation}\label{dispest4}
\norm{h(x-tv)}{L^{p}(\R^{d}_{x} ; L^{q} (V))}
\leq \frac{C(|V|)}{|t|^{d\left(\frac{1}{q}-\frac{1}{p}\right)}}
\norm{h(x)}{L^{q}(\R^{d}_{x})} .
\end{equation}

Next we recall the Strichartz estimates of \cite{CP}. 
\begin{proposition}\label{strichartz} (Strichartz estimates, \cite{CP}) 
Let $d\geq 2$ and let  $r,p,q,a \in [1,\infty]$ satisfy the conditions
\begin{equation}\label{conditions}
p\geq q, \quad \frac 2 r  = d \left(\frac{1}{q}-\frac{1}{p}\right)<1, \quad a = HM(p,q)\leq 2,
\end{equation}
where $HM$ denotes the harmonic mean.
\begin{enumerate}
\item
If $f(t,x,v)$ solves
\begin{equation}\label{NH}
\partial_{t} f + v\cdot\nabla_{x} f = g \ \ ,\ \  f(0,x,v)=0 ,
\end{equation}
then 
\begin{equation}\label{str1}
\norm{f}{L^{r}_{t} L^{p}_{x} L^{q}_{v} }\leq C \norm{g}{L^{r'}_{t} L^{q}_{x} L^{p}_{v}} .
\end{equation}

\item
If $f(t,x,v)$ solves
\begin{equation}\label{H}
\partial_{t} f + v\cdot\nabla_{x} f = 0 \ \ ,\ \  f(0,x,v)= f_{0}(x,v) ,
\end{equation}
then 
\begin{equation}\label{str2}
\norm{f}{L^{r}_{t} L^{p}_{x} L^{q}_{v} }\leq C \norm{f_{0}}{L^{a}_{x,v}}  .
\end{equation}
\end{enumerate}
\end{proposition}

\section{Global existence for arbitrarily large data}
In this Section we prove Theorems \ref{thm2} and \ref{thm1}. We start with  Theorem \ref{thm2}. Using the dispersion estimate gives rise to two norms, $\norm{\partial^\alpha S}{L^{p}}$
and $\norm{\rho}{L^{q}}$, see \eqref{secder3}. Of course each of them could be estimated in terms of $f$, but this
would result  in a quadratic term and would make the use of Gronwall's inequality impossible. However, thanks to
conservation of mass, we can choose $q=1$ (this  corresponds to velocity averaging) 
and bound $\norm{\rho}{L^{q}}$ a-priori. The norm   $\norm{\partial^\alpha S}{L^{p}}$ is then estimated 
using the well-known Calder\'on-Zygmund inequality if $|\alpha| = 2$, see \eqref{CZ} (see also Remark \ref{rem} at the end of this Section) and Young inequality if $|\alpha|\leq 1$.

We  use the standard abbreviations for mixed  spaces, 
for example $L^{p}_{x}L^{q}_{v}$ stands for  $L^{p}(\R^{d}_{x};L^{q}_{v}(V))$. 
In all cases $x$ varies in the whole space $\R^d$ while
$v$ and $v'$ are restricted in the bounded velocity space $V$.

\begin{proof}[Proof of Theorem \ref{thm2}.]
Fix $p$ and $q$ with $1\leq q \leq p \leq \infty$. 
Arguing as in \cite{CMPS} we have 
\begin{align}
f(t,x,v) & \leq f_{0}(x-tv,v) + C \int_{0}^{t} \rho(t-s,x-sv) dv  \nonumber \\ 
&  + C \sum_{|\alpha|\leq 2}
 \int_{0}^{t}\abs{ \partial^{\alpha} S(t-s,x-sv+v)} \rho( t-s,x-sv) ds \label{secder2}
\end{align}
therefore, using the dispersion estimate \eqref{dispest2}, we have
\begin{align}
& \norm{f(t,x,v)}{L^{p}_{x}L^{q}_{v}}  \leq  
 \norm{f_{0}(x-tv,v)}{L^{p}_{x}L^{q}_{v}} + C(|V|) 
  \int_{0}^{t} \frac{1}{s^{d\left(\frac{1}{q}-\frac{1}{p}\right)}} \norm{\rho(t-s,\cdot)}{L^{q}} ds  \nonumber \\
& + C \sum_{|\alpha|\leq  2}
\int_{0}^{t} \frac{1}{s^{d\left(\frac{1}{q}-\frac{1}{p}\right)}}
\norm{ \partial^{\alpha} S(t-s,x+v)  \rho( t-s,x)}{L^{q}_{x}L^{p}_{v}}  ds\nonumber \\
&\leq C_{0}(t) + C(|V|) \int_{0}^{t} \frac{1}{s^{d\left(\frac{1}{q}-\frac{1}{p}\right)}} \norm{\rho(t-s,\cdot)}{L^{q}} ds 
\nonumber\\
& + C \sum_{|\alpha|\leq 2}
\int_{0}^{t} \frac{1}{s^{d\left(\frac{1}{q}-\frac{1}{p}\right)}}
\norm{\partial^{\alpha} S(t-s, \cdot)}{L^{p}} \norm{\rho(t-s,\cdot)}{L^q} ds \label{secder3}
\end{align}
where we have set $C_{0}(t)=\norm{f_{0}(x-tv,v)}{L^{p}_{x}L^{q}_{v}}$.
Choose $q=1$ and  $p\in (1,\frac{d}{d-1})$.
Then by conservation of mass $\norm{\rho(t-s,\cdot)}{L^q} =M$. 
Using Young's inequality and conservation of mass for the derivatives of order one, we have 
\begin{equation}\label{gradient estimate}
\norm{\nabla S(t-s,\cdot)}{L^p} = C \norm{\rho(t-s,\cdot)*\nabla G}{L^p}\leq \norm{\rho(t-s,\cdot)}{L^1}
\norm{\nabla G}{L^p}=C M 
\end{equation}
where $G(x)=\frac{1}{4\pi}\int_{0}^{\infty} e^{-\pi\frac{|x|^2}{4s}-\frac{s}{4\pi}} s^{\frac{-d+2}{2}} \frac{ds}{s}$
is the Bessel potential,
and we get a similar estimate for $S$.
For the derivatives of order two we have (\cite{Stein}, p. 59, Proposition 3) 
\begin{multline}\label{CZ}
\norm{\partial_{ij} S(t-s)}{L^{p}} 
\leq C(d,p) \norm{\Delta S(t-s)}{L^p} \leq C(d,p) \norm{\rho(t-s)}{L^p} + C \norm{S(t-s)}{L^{p}}\\
\leq C(d,p) \norm{\rho(t-s)}{L^p} + CM.
\end{multline}
Therefore \eqref{secder3} gives
\[
\norm{\rho(t)}{L^{p}}
\leq C_{1}(t)  + C(d,p,M) \int_{0}^{t} \frac{1}{s^{d/p'}} \norm{\rho(t-s)}{L^p} ds .
\]
Since $d/p'<1$, we can use Gronwall's inequality to get
\begin{equation*}
\norm{\rho(t)}{L^{p}} \leq C(d,p,t,f_{0}) .
\end{equation*}
This completes the a-priori estimates. See remark \ref{averaging lemma}. 
\end{proof}

\begin{remark} We have chosen $\beta=1$ so that our $S$ decays sufficiently fast in order to apply the Calder\'on-Zygmund inequality. If $\beta=0$ we have $S=S^{s}+S^{l}\in L^p + L^\infty$ and we have no decay for $S^l$.
\end{remark}

The proof of Theorem \ref{thm1} 
uses the dispersion estimate of Proposition \ref{disp} as well as Young's convolution
inequality and the Hardy-Littlewood-Sobolev inequality. The dispersion estimate is used to
handle functions of $x-sv$ and $v$ which arise when we integrate the kinetic equation \eqref{Alt1},
see \eqref{f} and \eqref{f123} below. Each term in the right hand side of 
hypothesis \eqref{hyp1} requires an estimate in $L^{p}_{x}L^{q}_{v}$ for a certain range
of $p$ and $q$. Terms involving $x+v$ usually require small $p$  while terms involving
$x-v'$ require large $p$. The main difficulty then is to find one set of parameters
that makes both estimates work. To deal with this we will view the term 
$ \nabla S  (t-s, x-sv +v)  \rho(t-s,x-sv)$ as $\nabla S (t-s, x-(s-1)v) \rho(t-s,x-(s-1)v-v)$.
This shifting of the singularity from $s=0$ to $s=1$ (see \eqref{shift}) results in a redistribution
of  norms that allows us to estimate the terms involving $\nabla S (x+v)$
and $S(x+v)$ without any restrictions on the parameter $p$, and it  creates enough freedom 
so that, when we come to the more complicated estimates for $S(x-v')$, we
are able to find a pair $(p,q)$ that works for both.

\begin{proof}[Proof of Theorem \ref{thm1}] We shall only present a-priori estimates via a bootstrap argument 
for the solution $f$ of \eqref{Alt} in the space $ L^{p}(\R^{3}_{x};L^{q}_{v}(V))$.
The existence part of Theorem \ref{thm1} then follows by well-known methods, see Remark \ref{averaging lemma}.  
We present the proof in the more difficult case $\beta=0$.

Observe  that $S=S^{s} + S^{l}$ where $S^{s}(t)= \frac{1}{4\pi} \rho(t) * \frac{\mathbbm{1}_{|x|\leq 1}}{|x|}$
and $S^{l}(t)= \frac{1}{4\pi} \rho(t) * \frac{\mathbbm{1}_{|x|\geq 1}}{|x|}$. The long part $S^{l}(t)$ is a-priori bounded
thanks to conservation of mass:
\[
\abs{S^{l}(t,x)} \leq C \norm{\rho(t)}{L^{1}} 
\norm{\frac{\mathbbm{1}_{|x|\geq 1}}{|x|}}{L^{\infty}} \leq C M.
\] 
Similarly we split $\nabla S$ as $\nabla S  = \left(\nabla S\right)^{s} + \left(\nabla S\right)^{l}$
where $\left(\nabla S\right)^{s}(t)= \frac{1}{4\pi} \rho(t) * \frac{\mathbbm{1}_{|x|\leq 1}}{|x|^2}$
and $\left(\nabla S \right)^{l}(t)= \frac{1}{4\pi} \rho(t) * \frac{\mathbbm{1}_{|x|\geq 1}}{|x|^2}$ and show that 
$\left(\nabla S \right)^{l}$ is a-priori bounded. It follows that we may replace hypothesis \eqref{hyp1} by 
\begin{equation}\label{hyp1'}
0\leq T[S](t,x,v,v') \leq C\Big( 1+  S^{s}(t,x+v) + S(t,x-v') + \abs{ \left(\nabla S\right)^{s} (t,x+v)}\Big)
\end{equation}
where the new constant $C$ depends on the mass $M$.
For technical reasons it is more convenient not to split $S(t,x-v')$. Following the reasoning in 
\cite{CMPS} we estimate $f$ as follows:
\begin{equation}\label{f}
f(t,x,v) \leq C f_{0}(x-tv,v) + C \int_{0}^{t} \rho(t-s,x-sv) ds   + C \sum_{j=1}^{3} f_{j}(t,x,v)
\end{equation}
where
\begin{subequations}\label{f123}
\begin{align}
f_{1}(t,x,v)&= \int_{0}^{t}\int_{V} S^{s}(t-s, x-sv +v) f(t-s,x-sv,v') dv' ds \nonumber\\
& = \int_{0}^{t} S^{s}(t-s, x-sv +v) \rho(t-s,x-sv)  ds \label{f1}\\
f_{2}(t,x,v)&= \int_{0}^{t}\int_{V} S(t-s, x-sv -v') f(t-s,x-sv,v') dv' ds  \label{f2}\\
f_{3}(t,x,v)&= 
\int_{0}^{t} \abs{\left(\nabla S\right)^{s}(t-s, x-sv +v)} \rho(t-s,x-sv)  ds  \label{f3}
\end{align}
\end{subequations}

Fix $p$ and $q$ with $q\in [1, 3/2)$  such that 
\begin{equation}\label{integrability}
\lambda:=3 \left(\frac{1}{q}-\frac{1}{p}\right) < 1 , \ 1\leq q \leq p \leq \infty
\end{equation}
These restrictions on $p$ and $q$ will be enough for all estimates involving $f_{1}$ and $f_{3}$ 
but more restrictions 
will be imposed later when we estimate $f_{2}$ and we 
will want to know that there is a pair $(p,q)$ that satisfies all of them (\footnote{One choice of parameters
that works is: $p=\frac{9}{5}$, $q=\frac{9}{7}$.}).  
We start our estimates with  $f_{3}$. We have 
\begin{align*}
\norm{f_{3}(t,x,v)}{L^{p}_{x}L^{q}_{v}}& \leq \int_{0}^{t} 
\norm{\,  \left(\nabla S \right)^{s} (t-s, x-sv +v) \, \rho(t-s,x-sv)}{L^{p}_{x}L^{q}_{v}}  ds \\
&  \hspace{-0.3in} = \int_{0}^{t} 
\norm{\, \left(\nabla S \right)^{s}(t-s, x-(s-1)v) \, \rho(t-s,x-(s-1)v-v)}{L^{p}_{x}L^{q}_{v}}  ds .
\end{align*}
For fixed $t\geq 0$ and $s\in (0,t)$ we use the dispersion estimate \eqref{dispest3} with $t$ replaced by $s-1$ and
$h(x,v)=\abs{\left(\nabla S \right)^{s}(t-s, x)}  \rho(t-s,x-v)$ to get
\begin{align}
\norm{f_{3}(t,x,v)}{L^{p}_{x}L^{q}_{v}}&   \leq \int_{0}^{t} 
\frac{1}{|s-1|^{\lambda}} \norm{\,  \left(\nabla S \right)^{s}(t-s, x) \, 
\rho(t-s,x-v)}{L^{q}_{x}L^{p}_{v}} ds \label{shift} \\
&\leq \int_{0}^{t}\frac{1}{|s-1|^{\lambda}} 
\norm{\left(\nabla S \right)^{s} (t-s, \cdot)}{L^{q}}   \norm{\rho(t-s,\cdot)}{L^{p}}        ds .
\end{align}
Because $q<3/2$ the quantity $\norm{\left(\nabla S \right)^{s} (t-s, \cdot)}{L^{q}}$ 
is uniformly bounded.
Indeed, using Young's inequality we have 
\begin{equation}\label{f3B}
\norm{\left(\nabla S \right)^{s} (t-s, \cdot)}{L^{q}} \leq C \norm{\rho(t-s,\cdot)}{L^{1}}
\norm{\frac{\mathbbm{1}_{|x|\leq 1}}{|x|^2}}{L^{q}(\R^3)} \leq C(q) M . 
\end{equation}
On the other hand, since the velocity space is bounded, we have 
\[\norm{\rho(t-s,\cdot)}{L^{p}}=\norm{f(t-s,x,v)}{L^{p}_{x}L^{1}_{v}}\leq 
C(|V|,q)\norm{f(t-s,x,v)}{L^{p}_{x}L^{q}_{v}} . \] 
We conclude that 
\begin{equation}\label{f3C}
\norm{f_{3}(t,x,v)}{L^{p}_{x}L^{q}_{v}} \leq C(|V|,q) M  \int_{0}^{t} \frac{1}{|s-1|^{\lambda} } 
\norm{f(t-s,x,v)}{L^{p}_{x}L^{q}_{v}} ds .
\end{equation}
The estimate for $f_{1}$ is almost exactly the same. The only difference is that in the a-priori estimate
\eqref{f3B}  the norm $\norm{\frac{\mathbbm{1}_{|x|\leq 1}}{|x|^2}}{L^{q}(\R^3)}$ is replaced by 
$\norm{\frac{\mathbbm{1}_{|x|\leq 1}}{|x|}}{L^{q}(\R^3)}$ which is again finite because $q<3/2 < 3$. We get:
\begin{equation}\label{f1A}
\norm{f_{1}(t,x,v)}{L^{p}_{x}L^{q}_{v}} \leq C(|V|,q) M  \int_{0}^{t} \frac{1}{|s-1|^{\lambda} } 
\norm{f(t-s,x,v)}{L^{p}_{x}L^{q}_{v}} ds .
\end{equation}

\begin{remark}
Splitting in addition between small and long times ($s\gtrless 1/2$ for example) we end up 
with a priori estimates without any restriction on the exponent $p$. But this technical improvement is not relevant in this proof.
\end{remark}

Next we estimate $f_{2}$. We start with some numerology which we explain later. 
Fix $q\in (1,3/2)$. There exists a $p\in (3/2, 3)$ such that 
\begin{equation}\label{numerology1} 3 \left( \frac{1}{q} - \frac{1}{p}\right) = 1 - \frac{3p'}{(q')^2} .\end{equation}
To see this write \eqref{numerology1} as 
$\delta(p):=1 - \frac{3p'}{(q')^2} - 3 \left( \frac{1}{q} - \frac{1}{p}\right)=0$
and think of this expression as a continuous function of the variable 
$p\in [ 3/2 , 3 ]$. For $p=3/2$ we have
$\delta(3/2)=1-\frac{9}{(q')^2} - 3 \left( \frac{1}{q} - \frac{2}{3}\right)= \frac{3(q'-3)}{(q')^2}>0$. 
On the other hand for $p=3$ we have 
$\delta(3)=1 - \frac{9}{2(q')^2} - 3 \left( \frac{1}{q} - \frac{1}{3}\right)<- \frac{9}{2(q')^2}<0$. 
The existence of $p$ follows. Notice that with this
choice of $p$ and $q$ we still have $1\leq q \leq p \leq \infty$ and moreover the integrability condition 
$\lambda<1$  is satisfied thanks to \eqref{numerology1}.

\begin{remark} \label{rem:Admissible exponents}
In fact we are allowed to choose $q$ and $p$ to be in the following set of exponents
\begin{equation} \label{numerology1bis}  \mathcal{A} = \Big\{ p'\geq 1,\ q'\geq 1 \ \Big| \ q'> p' \ , \  
3\left( \frac{1}{q} - \frac{1}{p}\right) + \frac{3p'}{(q')^2} \leq  1 
 \ , \ \frac1{q} - \frac1{p}<\frac1{3} \Big\}.\end{equation}
This set of admissible exponents for the estimate of $\|f_2(t,x,v)\|_{L^p_x L^q_v}$ is plotted in figure \ref{Admissible exponents}a in the coordinates $(q',p')$. The key point is that it intersects the constrain $\{q'> 3\}$ which comes from the estimates on $f_1$ and $f_{3}$. This can be done only by decoupling $p$ and $q$. \\
Assuming some linear contribution of $\nabla S(x-v')$ in the turning kernel bound \eqref{hyp1} would have lead to the set represented in figure \ref{Admissible exponents}b. The latter does not intersect the half-plane $\{q'>3\}$.
\begin{figure}
\begin{center}
\includegraphics[width = .45\linewidth]{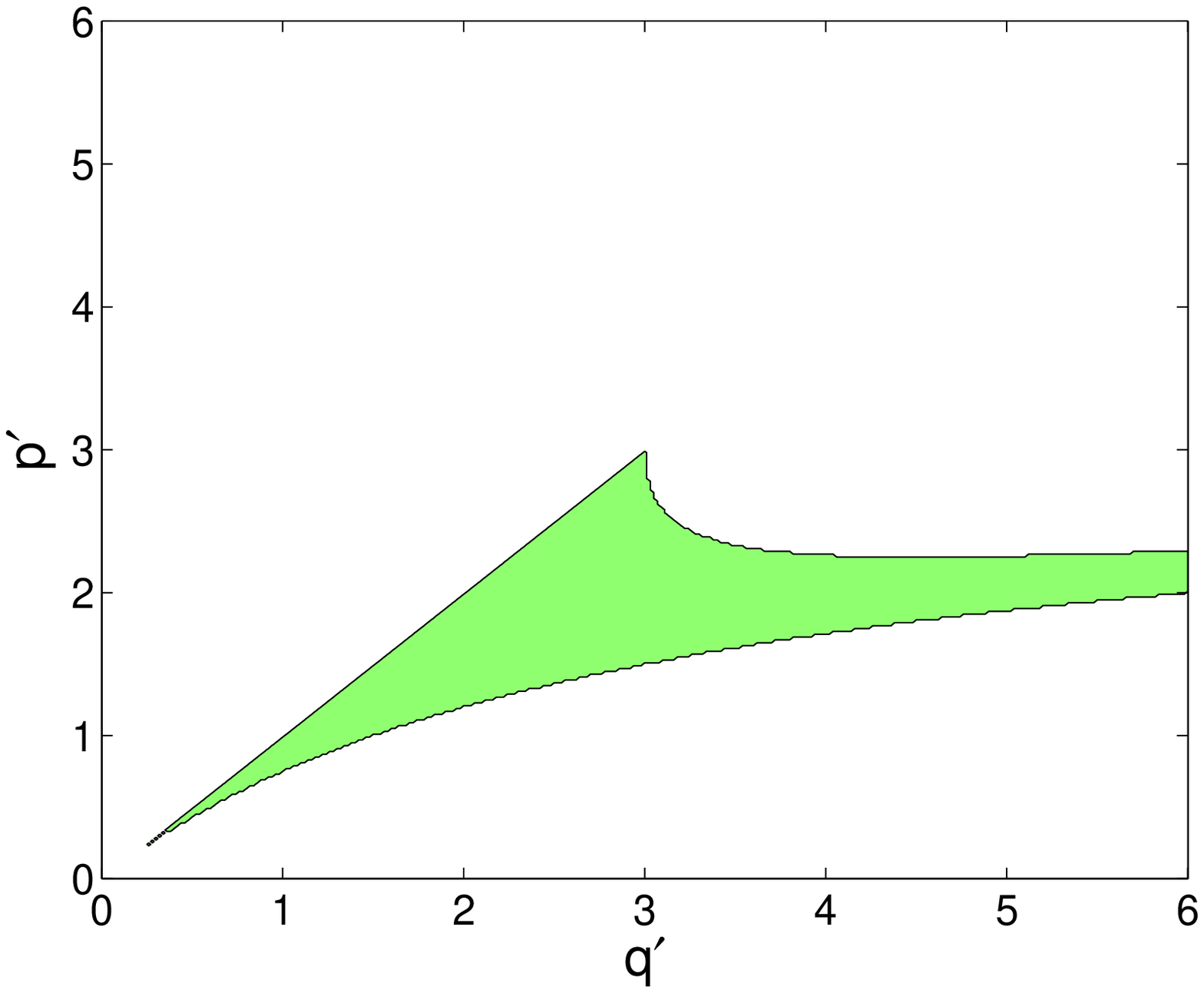}(a)\
\includegraphics[width = .45\linewidth]{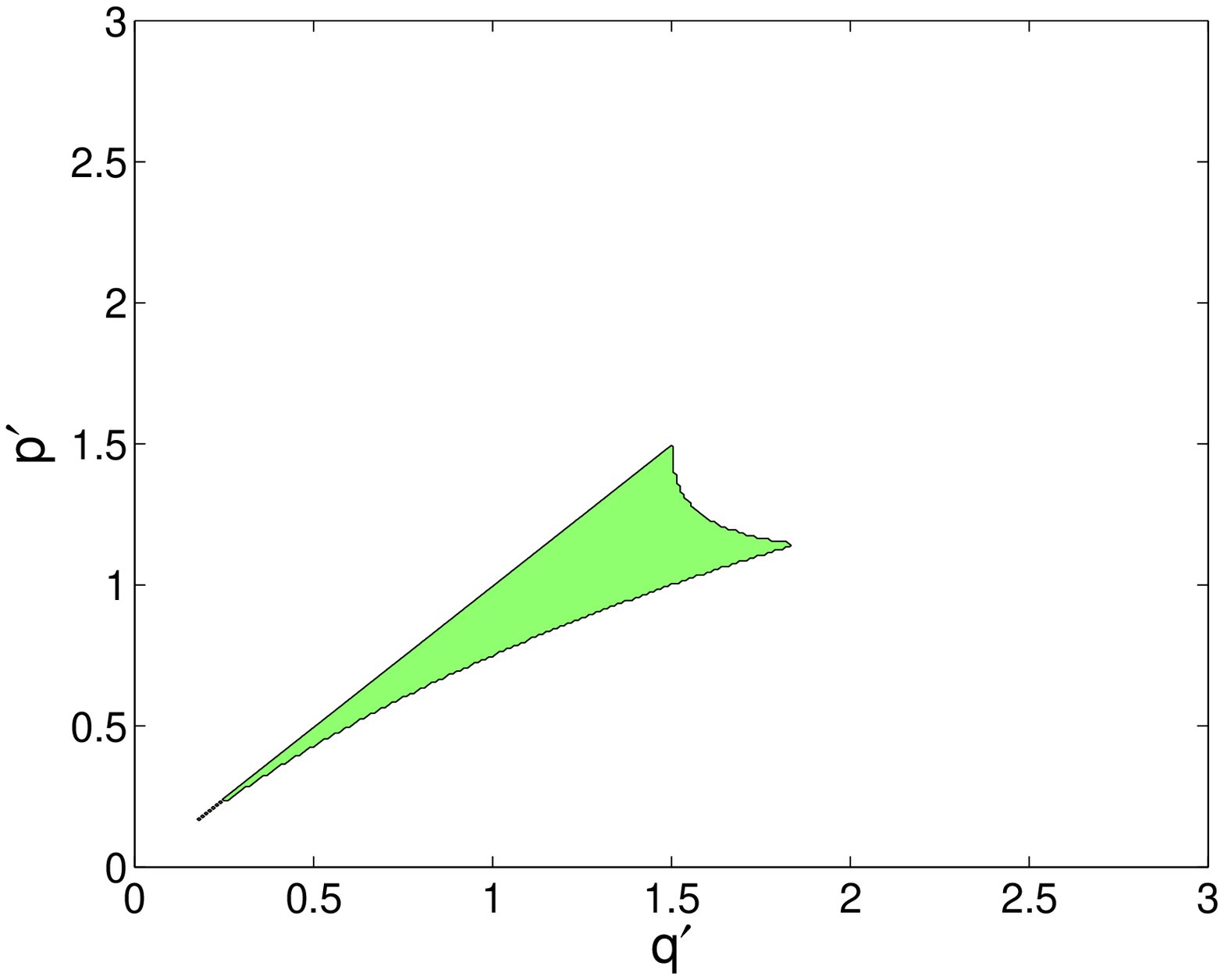}(b) 
\caption{Set of admissible exponents $(q',p')$ for the estimate of $\|f_2(t,x,v)\|_{L^p_x L^q_v}$,
corresponding to $S(x-v')$ (a) and $\nabla S(x-v')$ (b) respectively.}
\label{Admissible exponents}
\end{center}
\end{figure}

\end{remark}

Finally we define $\theta\in (0,1)$,  $c\in (1,q)$ and $b\in (1,c')$ by the following relations:
\begin{equation}\label{numerology2}
\frac{1}{q}=1-\theta + \frac{\theta}{p}\ ,\ \frac{1}{c}=1-\theta + \frac{\theta}{q}\ ,\ 
\frac{1}{b}+\frac{1}{c}=\frac{5}{3} .
\end{equation} 
Proceeding with the estimates we have 
\begin{equation}\label{f2A}
\norm{f_{2}(t,x,v)}{L^{p}_{x}L^{q}_{v}} \leq \int_{0}^{t} 
\norm{\, \int_{V}  S (t-s, x-sv -v') \, f(t-s,x-sv,v')\, dv'}{L^{p}_{x}L^{q}_{v}}  ds .
\end{equation}
For fixed $t\geq 0$ and $s\in (0,t)$ use the dispersion estimate \eqref{dispest4} with 
\[h(x)=\int_{V} S(t-s,x-v') f(t-s,x, v') dv'\]
 to get
\begin{multline}\label{f2B}
\norm{\, \int_{V}  S (t-s, x-sv -v') \, f(t-s,x-sv,v')\, dv'}{L^{p}_{x}L^{q}_{v}}  \\
  \leq \frac{1}{s^\lambda} \norm{\int_{V} S(t-s,x-v') f(t-s,x, v') dv'}{L^{q}_{x}} .
\end{multline}
By H\"older's inequality
\begin{equation}\label{f2C}
\int_{V} S(t-s,x-v') f(t-s,x, v') dv'  \leq 
\norm{S(t-s,\cdot)}{L_x^{c'}} \norm{f(t-s,x,v')}{L^{c}_{v'}} .
\end{equation}
Since $1<b<c'<\infty$ and $\frac{1}{b}-\frac{1}{c'}=\frac{2}{3}$ we can apply 
the Hardy-Littlewood-Sobolev inequality (see for instance \cite{Stein}, Theorem 1, page 199) to get
\begin{equation}\label{f2D}
\norm{S(t-s,\cdot)}{L^{c'}}= C \norm{\rho(t-s,\cdot) * \frac{1}{|x|}}{L^{c'}(\R^3)}
 \leq C \norm{\rho(t-s,\cdot)}{L^{b}} .
\end{equation}
It is easy to check that $1<b<p$ 
(\footnote{$\frac{1}{b}-\frac{1}{p}=\frac{5}{3}-\frac{1}{c}-\frac{1}{p}
=\frac{5}{3} -1 + \theta - \frac{\theta}{q} -\frac{1}{p}=
\left(\frac{2}{3}-\frac{1}{p}\right) + \frac{\theta}{q'} >0$.}), 
therefore if we define $\epsilon \in (0,1)$ by 
\begin{equation}\label{numerology3}
\frac{1}{b}=1-\epsilon + \frac{\epsilon}{p}
\end{equation}
we can use interpolation and  conservation of mass to obtain 
\begin{equation*}
\norm{\rho(t-s)}{L^{b}}\leq 
\norm{\rho(t-s)}{L^{1}}^{1-\epsilon}
\norm{\rho(t-s)}{L^{p}}^{\epsilon} 
\leq M^{1-\epsilon} \norm{\rho(t-s)}{L^{p}}^{\epsilon} .
\end{equation*}
We have shown that
\begin{equation*}
 \int_{V} S(t-s,x-v') f(t-s,x, v') dv'  \leq C M^{1-\epsilon} \norm{\rho(t-s,\cdot)}{L^{p}}^{\epsilon}
\norm{f(t-s,x,v')}{L^{c}_{v'}}\ ,
\end{equation*}
and as a product we obtain
\begin{multline}\label{f2G}
\norm{ \int_{V} S(t-s,x-v') f(t-s,x, v') dv' }{L^{q}_{x}} \\
 \leq C M^{1-\epsilon} \norm{\rho(t-s,\cdot)}{L^{p}}^{\epsilon}
\norm{f(t-s,x,v)}{L^{q}_{x}L^{c}_{v}} .
\end{multline}
We aim to interpolate the  $L^{q}_{x}L^{c}_{v}$-norm between $L^1L^1$ and $L^p L^q$ in order to 
conclude with a Gronwall lemma. This is achieved thanks
to the first two relations in \eqref{numerology2}.
\begin{align}
\norm{f(t-s,x,v)}{L^{q}_{x}L^{c}_{v}}&  \leq \norm{f(t-s,x,v)}{L^{1}_{x}L^{1}_{v}}^{1-\theta}
\norm{f(t-s,x,v)}{L^{p}_{x}L^{q}_{v}}^{\theta} \nonumber \\
& \leq M^{1-\theta} \norm{f(t-s,x,v)}{L^{p}_{x}L^{q}_{v}}^{\theta}  .  \label{f2H}
\end{align}
Using this estimate together with 
$\norm{\rho(t-s,\cdot)}{L^{p}} \leq C(|V|,q)\norm{f(t-s,x,v)}{L^{p}_{x}L^{q}_{v}}$ 
into \eqref{f2G}, we get
\begin{equation}\label{f2I}
\norm{ \int_{V} S(t-s,x-v') f(t-s,x, v') dv' }{L^{q}_{x}} \leq C(|V|,q) M^{2-(\epsilon + \theta)}
\norm{f(t-s,x,v)}{L^{p}_{x}L^{q}_{v}}^{\epsilon + \theta} .
\end{equation}
We can now argue that we opted for \eqref{numerology1} to ensure that 
$\epsilon + \theta =1\,$
(\footnote{ We have $\theta=\frac{p'}{q'}$ and $\epsilon=\frac{p'}{b'}$ therefore 
$\epsilon + \theta =1$ is equivalent to $\frac{1}{q'}+\frac{1}{b'}=\frac{1}{p'}$. We calculate
$\frac{1}{p'} - \frac{1}{q'} - \frac{1}{b'} = \left(\frac{1}{q}-\frac{1}{p}\right) -1 + \frac{1}{b}$.
Now $-1+\frac{1}{b}=-1+\frac{5}{3}-\frac{1}{c}=\frac{2}{3}-\left(1-\theta + \frac{\theta}{q}\right)=
-\frac{1}{3} + \frac{\theta}{q'} =-\frac{1}{3} + \frac{p'}{(q')^2}$, 
therefore 
$\frac{1}{p'} - \frac{1}{q'} - \frac{1}{b'}=\left(\frac{1}{q}-\frac{1}{p}\right) -\frac{1}{3} + \frac{p'}{(q')^2}$
which is equal to zero thanks to \eqref{numerology1}. }). 
Therefore
\begin{equation}\label{f2J}
\norm{ \int_{V} S(t-s,x-v') f(t-s,x, v') dv' }{L^{q}_{x}} \leq C(|V|,q) M
\norm{f(t-s,x,v)}{L^{p}_{x}L^{q}_{v}} .
\end{equation}
From \eqref{f2A}, \eqref{f2B} and \eqref{f2J}, we conclude that
\begin{equation}\label{f2K}
\norm{f_{2}(t,x,v)}{L^{p}_{x}L^{q}_{v}} \leq C(|V|,q) M \int_{0}^{t} 
\frac{1}{s^\lambda} \norm{f(t-s,x,v)}{L^{p}_{x}L^{q}_{v}}  ds .
\end{equation}
This completes the estimate of $f_{2}$. Finally we have to estimate the first two terms in the right hand side of 
\eqref{f}. For the first term we have by our hypothesis on the initial data that 
$\norm{f_{0}(x-tv,v)}{L^{p}_{x}L^{q}_{v}} =:C_{0}(t)<\infty$ for all $t$. 
We can use dispersion and interpolation for the second term, leading to
\[ 
\norm{\int_{0}^{t} \rho(t-s, x-sv) ds}{L^{p}_{x}L^{q}_{v}} \leq 
\int_{0}^{t} s^{-\lambda}  \norm{\rho(t-s) }{L^{q}_x} ds \leq C \int_{0}^{t} s^{-\lambda}  \norm{\rho(t-s) }{L^{p}_x}^{\theta} ds,
\]
where $\theta$ has already been defined in \eqref{numerology2}.


Putting everything together we conclude that 
\begin{multline}\label{GronwallA}
\norm{f(t,x,v)}{L^{p}_{x}L^{q}_{v}} \leq C_{0}(t) + C(|V|,q) M \int_{0}^{t} 
K(s) \norm{f(t-s,x,v)}{L^{p}_{x}L^{q}_{v}}  ds \\  
+C \int_{0}^{t} s^{-\lambda}  \norm{f(t-s) }{L^{p}_xL^q_v}^{\theta} ds,
\end{multline}
where $K(s)=  1 + \frac{1}{s^\lambda} + \frac{1}{|s-1|^{\lambda}}$. Since $\lambda < 1$
we can apply Gronwall's inequality to obtain 
\begin{equation}\label{GronwallB}
\norm{f(t,x,v)}{L^{p}_{x}L^{q}_{v}} \leq C(|V|,q,t,f_{0}) .
\end{equation}

\end{proof}

\begin{remark}\label{rem} It would be interesting to know whether, in the case $\beta=1$, the hypotheses \eqref{hyp1} and \eqref{hyp2} can be 
combined into the single assumption: 
\begin{multline*}
0\leq T[S](t,x,v,v') \leq C \left(  1+   S(t,x+v) + S(t,x-v') + \abs{\nabla S (t,x+v)}\right)\\
 + C \sum_{i,j=1}^{3}\abs{\partial_{ij}S(t,x+v)}.
\end{multline*}
The obstruction in our estimates is that  the proof of Theorem \ref{thm2} requires $q=1$, so that the 
norm $\norm{\rho}{L^{q}}$
in \eqref{secder3} can be estimated a-priori, while the estimates for $f_{2}$ in the proof of Theorem \ref{thm1}
do not work with $q=1$ because it gives $b=3/2$, $c=1$ which is not allowed in the HLS inequality in \eqref{f2D}.
\end{remark}

\bigskip

\section{Global existence for small data in the critical norm}
Strichartz estimates have been very successful in dealing with many classes of nonlinear Schr\"odinger, wave and 
other dispersive equations. Typically they are used to show either local existence of solutions with
low regularity data or global existence under an additional smallness assumption on the initial data,
see \cite{Tao}.

\begin{proof}[Proof of Theorem \ref{thm3}]
To simplify the notation we use again the standard abbreviations for mixed spaces, for example
$L^{r}_{t}L^{p}_{x}L^{q}_{v}$. In all cases the space variable $x$ runs through all of $\R^3$
and the velocity variables $v$ and $v'$  always vary in the velocity space $V$.

Observe first that hypothesis \eqref{hyp3} implies that for all  $p_{1},p_{2},p_{3} \in [1,\infty]$
 with $ p_{1} \geq p_{2} , p_{3} $ we have
\begin{equation} \label{hyp4}
\norm{T[S](t,x,v,v')}{L^{p_{1}}_{x} L^{p_{2}}_{v} L^{p_{3}}_{v'}}\leq C(|V|, p_{2},p_{3})\left[
\norm{S(t,\cdot)}{L^{p_{1}}} + \norm{\nabla S(t,\cdot)}{L^{p_{1}}} \right] .
\end{equation} 
Indeed, since $p_{1}\geq p_{2} , p_{3}$, we can use Minkowski's inequality to obtain
\[\norm{S(t,x+v)}{L^{p_{1}}_{x}L^{p_{2}}_{v}L^{p_{3}}_{v'}}\leq 
\norm{S(t,x+v)}{L^{p_{2}}_{v}L^{p_{3}}_{v'}L^{p_{1}}_{x}} = 
C(|V|)\norm{S(t,\cdot)}{L^{p_{1}}}\] 
with similar estimates for all other terms in the right hand side of \eqref{hyp3}.
From now on the proof will use  estimate \eqref{hyp4} instead of hypothesis \eqref{hyp3}. 
We  present  a bootstrap argument for the solution $f$ in the space 
$L^{3}_{t}L^{p}_{x}L^{q}_{v}$. The existence result of Theorem \ref{thm3} then follows by 
standard methods. 

As usual we have:
\[f(t,x,v) \leq f_{1}(t,x,v) + f_{2}(t,x,v)\]
where $f_{1}(t,x,v)$ solves 
\[\partial_{t} f_{1} + v\cdot\nabla_{x} f_{1} = 0 \ \ ,\ \  f_{1}(0,x,v)= f_{0}(x,v)\]
and $f_{2}(t,x,v)$ solves
\[\partial_{t} f_{2} + v\cdot\nabla_{x} f_{2} = \int_{V} T[S] f' dv'\ \ ,\ \  f_{2}(0,x,v)= 0 .\] 
Recall that $a\in [3/2,2]$. Choose $r=3$ and define $p\in[ 9/5 , 18/7]$ and $q\in [ 9/7 , 18/11]$ by 
$\frac{1}{p}=\frac{1}{a}-\frac{1}{9}$ and $\frac{1}{q}=\frac{1}{a}+\frac{1}{9}$. 
It is easy to verify that the  quadruplet $(r,p,q,a)$ satisfies  the conditions \eqref{conditions} 
required for applying the Strichartz estimates. Apply  estimate \eqref{str1} to $f_1(t,x,v)$ 
and  estimate \eqref{str2} to $f_{2}(t,x,v)$ to get:
\begin{align}
\norm{f}{L^{3}_{t} L^{p}_{x} L^{q}_{v}} &  \leq \norm{f_1}{L^{3}_{t}L^{p}_{x}L^{q}_{v}} + 
\norm{f_2}{L^{3}_{t} L^{p}_{x} L^{q}_{v}} \\
& \leq C \norm{f_{0}}{L^{a}_{x,v}} + C \norm{\int_{V} T[S] f' dv' }{L^{3/2}_{t}L^{q}_{x}L^{p}_{v}} . 
\label{3str}
\end{align}
To estimate the last term in \eqref{3str}, apply first H\"older's inequality to get:
\[\int_{V} T[S](t,x,v,v') f(t,x,v') dv' \leq 
\norm{T[S](t,x,v,v')}{L^{q'}_{v'}} \norm{f(t,x,v')}{L^{q}_{v'}} .\]
Taking the $L^{p}_{v}$-norm of both sides we find:
\[\norm{ \int_{V} T[S] f' dv' }{L^{p}_{v}} 
\leq \norm{T[S](t,x,v,v')}{L^{p}_{v} L^{q'}_{v'}} \norm{f(t,x,v')}{L^{q}_{v'}} .\]
Taking next the  $L^{q}_{x}$-norm of both sides  and using H\"older's inequality  with 
$\frac{1}{q}=\frac{1}{p}+\frac{2}{9}$ we find:
\begin{equation}\label{3m14}
\norm{ \int_{V} T[S] f' dv'}{L^{q}_{x} L^{p}_{v}}\leq 
\norm{T[S](t,x,v,v')}{L^{9/2}_{x} L^{p}_{v} L^{q'}_{v'}}
\norm{f(t,x,v')}{L^{p}_{x} L^{q}_{v'}} .
\end{equation}
It is easy to check that $\frac{9}{2} \geq p$ and $\frac{9}{2} \geq q'$, therefore we can use 
\eqref{hyp4}  to get
\begin{align*}
\norm{T[S](t,x,v,v')}{L^{9/2}_{x}(L^{p}_{v}(L^{q'}_{v'}))} & \leq C(|V|,p,q')\left[ 
\norm{S(t,x)}{L^{9/2}_{x}} +  \norm{\nabla S(t,x)}{L^{9/2}_{x}} \right] \\
&= C(|V|,p,q')\left[ \norm{ G * \rho(t)}{L^{9/2}_{x}} + 
 \norm{\nabla G * \rho(t)}{L^{9/2}_{x}} \right]
\end{align*}
where $G$ is the Bessel potential (see the proof of theorem \ref{thm2}).

If $\frac{9}{5} < p$ we proceed by using Young's inequality. One 
can show that $G\in L^{b}$ for all $b<3$ and that $\nabla G \in L^{b}$ for all $b<\frac{3}{2}$.
Define $b$ by  $1+\frac{2}{9}=\frac{1}{b} + \frac{1}{p}$. Then $\frac{6}{5}\leq b<\frac{3}{2}$. 
We get
\begin{align*}
 \norm{G  * \rho}{L^{9/2}_{x}} + 
\norm{\nabla G  * \rho}{L^{9/2}_{x}} 
&\leq \norm{G}{L^b} \norm{\rho}{L^{p}_{x}} +\norm{\nabla G}{L^b} \norm{\rho}{L^{p}_{x}} \\
&\leq C(b)  \norm{\rho(t,x)}{L^{p}_{x}}\\
&\leq C(b,q,|V|) \norm{f(t,x,v)}{L^{p}_{x} L^{q}_{v}} .
\end{align*}
If $\frac{9}{5}=p$ we use Young's inequality for the $G$-term and the HLS inequality for the $\nabla G$-term.
Defining $b$ as above now gives $b=\frac{3}{2}<3$ therefore
\begin{equation}\label{3hls1}
\norm{G * \rho}{L^{9/2}_{x}}\leq \norm{G}{L^{3/2}}\norm{\rho}{L^{p}_{x}} \leq C  \norm{\rho}{L^{p}_{x}}
\leq C(q,|V|) \norm{f(t,x,v)}{L^{p}_{x} L^{q}_{v}} .
\end{equation}
One can show that $\abs{\nabla G(x)}\leq  \frac{C}{|x|^2}$ for all $x$. 
Therefore, by HLS,
\[
\norm{\nabla G * \rho}{L^{9/2}_{x}} 
\leq  \norm{\frac{C}{|x|^2} * \rho}{L^{9/2}_{x}}
 \leq C \norm{\rho}{L^{9/5}_{x}}
= C \norm{\rho}{L^{p}_{x}}
\leq C(q,|V|) \norm{f(t,x,v)}{L^{p}_{x} L^{q}_{v}} .
\]

\smallskip

The above argument shows that
\begin{equation}\label{3m15}
\norm{T[S](t,x,v,v')}{L^{9/2}_{x} L^{p}_{v} L^{q'}_{v'}} \leq C(a,|V|)  \norm{f(t,x,v)}{L^{p}_{x} L^{q}_{v}} .
\end{equation}
Using \eqref{3m15} into \eqref{3m14} we get:
\begin{equation}
\norm{\int_{V} T[S] f' dv'}{L^{q}_{x} L^{p}_{v}}\leq C(a,|V|)
\norm{f(t,x,v)}{L^{p}_{x}L^{q}_{v}}^2.
\end{equation}
Taking the  $L^{3/2}_{t}$-norm of both sides we obtain:
\begin{equation}
\norm{\int_{V} T[S] f' dv'}{L^{3/2}_{t} L^{q}_{x} L^{p}_{v}}\leq 
\norm{\norm{f(t,x,v)}{L^{p}_{x}(L^{q}_{v})}^2}{L^{3/2}_{t}}=\norm{f(t,x,v)}{L^{3}_{t} L^{p}_{x} L^{q}_{v}}^2.
\end{equation}
Using this in \eqref{3str} we find
\begin{equation}
\norm{f}{L^{3}_{t} L^{p}_{x}L^{q}_{v}}  \leq C \norm{f_{0}}{L^{a}(\R^6)}
+ C \norm{f(t,x,v)}{L^{3}_{t} L^{p}_{x} L^{q}_{v}}^2 .
\end{equation}

This completes the a-priori estimates  which enable to bootstrap for small initial data. 
See remark \ref{averaging lemma}. 
\end{proof}

\begin{remark}\label{numerology4}
The proof of Theorem \ref{thm3} works in $d=4$ dimensions too. One may choose for instance
$(q,p,r,a)=(3,12/5,12/7,2)$. Notice that $a=\frac{d}{2}$. It also works in dimension $d=2$, 
however, in this case a better result (global existence without a smallness assumption)
is available in \cite{HKS1}. \\
Using the same method we can prove local existence for large data.
\end{remark}

\begin{remark}\label{averaging lemma}
We have proved a priori estimates for the IVP \eqref{Alt}. We can prove the existence of weak solutions using regularization and compactness. In particular the compactness can be gained using averaging lemmas (see \cite{P2} for instance) provided we get some a priori estimate on the $L^pL^q$-norm of $f$ with $q>1$. This has been obtained in the proofs of theorems \ref{thm1}
and \ref{thm3}, whereas in theorem \ref{thm2} an additional bootstrap step is needed. Of course continuity of
$T[S]$ in spaces $L^p_{loc}$ is needed for passing to the limit in all cases.

\end{remark}

\section*{Conclusions and perspectives} In this paper we have considered a number of classes of turning kernels in the kinetic model of chemotaxis. We have proved global existence for arbitrarily large data using dispersion estimates for several of them, and, using Strichartz estimates, we have obtained global existence for small solutions in the most difficult case of a turning kernel that involves $\abs{\nabla S (t,x+v)}+\abs{\nabla S (t,x-v')}$ (Theorem 3). The  de-localization induced by $v$ or $v'$ in these formula is fundamental both for mathematical theory and biophysical interpretation. However, several questions remain that show that the present theory still needs to be improved. We would like to mention a few of them. 

At first, obviously is  the case of large initial data in Theorem \ref{thm3} which remains open. Notice that the time integrability in the Strichartz estimates implies some decay to zero at infinity which is only possible for small initial data, as we know from the Keller-Segel system \cite{BDP,CPZ}.

 A second question is to include some of these examples in a more general assumption such as
$$
\|T[S](t,x,v',v) \|_{L^\infty_{loc}} \leq C \|S \|_{L^\infty_{loc}}.
$$
Because this does not include directly de-localization, the methods used here do not apply for global existence. 

More related to  biophysical interpretation there is a  third question:
unlike in the Keller-Segel model -- where having turning kernels of the form $S(t,x+v')$ and $\nabla S(t,x+v')$
gives a repellent drift \cite{CMPS} and then there is no blow-up, and the existence theory is much simpler -- in the arguments carried out in the proof of the results we do not see why it should be better to have turning kernels of the latter form. 
\\
\\
 

\noindent {\bf Acknowledment.} The authors would like to acknowledge support from European Networks MRNTN-CT-2004-503661,  HPRN-CT-2002-00282, and HPRN-CT-2001-00273 - HARP.

\end{document}